\documentclass{article}

\usepackage[utf8]{inputenc}
\usepackage{CJKutf8}

\usepackage[inline]{asymptote}

\usepackage[version=4]{mhchem}

\usepackage{array}

\usepackage{graphicx}
\graphicspath{ {./images/} }

\usepackage[usenames,dvipsnames]{color}
\definecolor{dkgreen}{rgb}{0,0.6,0}
\definecolor{gray}{rgb}{0.5,0.5,0.5}
\definecolor{mauve}{rgb}{0.58,0,0.82}
\definecolor{lavender}{rgb}{0.8,0.56,0.98}
\definecolor{emerald}{RGB}{21, 89, 26}

\usepackage{tikz}
\usetikzlibrary{arrows.meta,shapes.misc,decorations.pathmorphing,calc,bending,decorations.markings,shapes.geometric}
\usepackage[makeroom]{cancel}

\usepackage[normalem]{ulem}

\usepackage{amsmath}
\usepackage{mathtools}
\usepackage{amssymb}
\usepackage{wasysym}
\usepackage{amsfonts}
\usepackage{mathrsfs}
\DeclareMathAlphabet{\mathpzc}{OT1}{pzc}{m}{it}

\usepackage{bigints}

\usepackage{relsize}

\usepackage[english]{babel} 
\usepackage[autostyle, english = american]{csquotes}
\MakeOuterQuote{"} % sets the double quotation mark as an outer quote

\usepackage{hyperref}
\hypersetup{
    colorlinks,
    linkcolor={red!50!black},
    citecolor={blue!50!black},
    urlcolor={blue!80!black}
}

\usepackage{geometry}
\usepackage{mdframed}
\usepackage{calc}
\definecolor{mytitlecolor}{RGB}{0,102,204}
\usepackage{tcolorbox}
\definecolor{exercolor}{RGB}{0,200,115}

\usepackage{tcolorbox}
\definecolor{examplecolor}{RGB}{161,161,161}

\usepackage{mdframed}
\usepackage{calc}
\definecolor{myDefinitionColor}{RGB}{255, 104, 0}
\usepackage{listings}
\usepackage{svg}
\makeatletter
\def\smallunderbrace#1{
  \mathop{
    \vtop{
      \m@th
      \ialign{##\crcr
        $\hfil\displaystyle{#1}\hfil$\crcr
        \noalign{\kern3\p@\nointerlineskip}
        \tiny\upbracefill\crcr
        \noalign{\kern3\p@}
      }
    }
  }\limits
}
\makeatother

\usepackage[d]{esvect}

\usepackage{pgf}
\newcommand*{\SectorRadius}{1ex}
\newcommand*{\SectorHalfAngle}{45}
\newcommand*{\SectorLineWidth}{.4pt}
\newcommand*{\sector}{%
  \begin{pgfpicture}
    \pgfpathmoveto{\pgforigin}%
    \pgfpathlineto{\pgfpointpolar{90-(\SectorHalfAngle)}{\SectorRadius}}%
    \pgfarc{90-(\SectorHalfAngle)}{90+\SectorHalfAngle}{\SectorRadius}%
    \pgfpathclose
    \pgfsetlinewidth{\SectorLineWidth}%
    \pgfusepath{stroke}%
  \end{pgfpicture}%
}

\newcommand{\semicirc}{\raisebox{-0.13em}{\rotatebox[origin=c]{90}{$\Rightcircle$}}\,}

\usepackage{biblatex}
\DeclareCiteCommand{\supercite}[\mkbibsuperscript]
  {\iffieldundef{prenote}
     {}
     {\BibliographyWarning{Ignoring prenote argument}}%
   \iffieldundef{postnote}
     {}
     {\BibliographyWarning{Ignoring postnote argument}}}
  {\usebibmacro{citeindex}%
   \bibopenbracket\usebibmacro{cite}\bibclosebracket}
  {\supercitedelim}
  {}
\let\cite=\supercite

\usepackage{etoolbox}
\makeatletter
\patchcmd{\l@section}
    {\hfil}
    {\leaders\hbox{\normalfont$\hbox{.}\mkern \@dotsep mu$}\hfill}
    {}{}
\makeatother
\usepackage{multicol}
\begin{document}

% Title
\begin{center}\vspace{-1cm}
    \textbf{\LARGE Closed Form for Half-Area Overlap Offset of 2 Unit Disks}\\
    Max Chicky Fang 
    ~
\end{center}

\begin{center}
\subsection*{Abstract}

    The separation between the centers of two unit circles such that their overlapping area is exactly half of each's area is known to be around $\mathpzc{D}_{\text{DHA}} = 0.8079455\dots$ (OEIS \href{https://oeis.org/A133741}{A133741}). However, no closed form of $\mathpzc{D}_{\text{DHA}}$ is known. Here, we determine its closed form representation in terms of the inverse regularized beta function.
\end{center}

\section{Introduction} 

\subsection{Setup}

    Consider the diagram below.
    
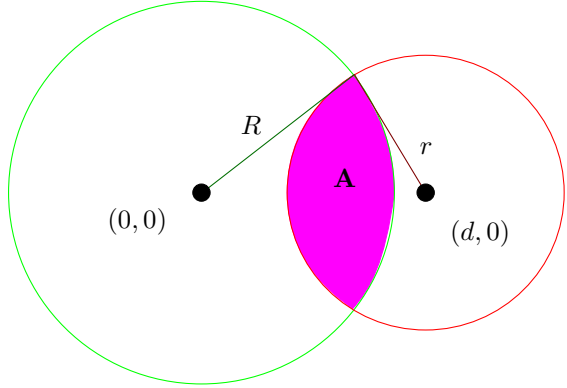
\begin{center}
    \begin{tikzpicture}[x=0.75pt,y=0.75pt,yscale=-1,xscale=1]
        %uncomment if require: \path (0,300); %set diagram left start at 0, and has height of 300
        
        %Shape: Polygon Curved [id:ds07395154754330058] 
        \draw  [color={rgb, 255:red, 255; green, 0; blue, 253 }  ,draw opacity=1 ][fill={rgb, 255:red, 255; green, 0; blue, 253 }  ,fill opacity=1 ] (212.28,145.7) .. controls (212.8,121.8) and (221.8,102.8) .. (245.81,86.27) .. controls (258.8,109.8) and (264.8,117.8) .. (265.77,145.7) .. controls (260.8,181.8) and (253.8,193.8) .. (245,204) .. controls (231.8,194.8) and (212.8,177.8) .. (212.28,145.7) -- cycle ;
        %Shape: Ellipse [id:dp3499342112834707] 
        \draw  [color={rgb, 255:red, 0; green, 255; blue, 0 }  ,draw opacity=1 ] (73.27,145.7) .. controls (73.27,92.54) and (116.36,49.45) .. (169.52,49.45) .. controls (222.68,49.45) and (265.77,92.54) .. (265.77,145.7) .. controls (265.77,198.85) and (222.68,241.94) .. (169.52,241.94) .. controls (116.36,241.94) and (73.27,198.85) .. (73.27,145.7) -- cycle ;
        %Shape: Ellipse [id:dp8007078804190029] 
        \draw  [color={rgb, 255:red, 255; green, 0; blue, 0 }  ,draw opacity=1 ] (212.28,145.7) .. controls (212.28,107.53) and (243.21,76.59) .. (281.38,76.59) .. controls (319.54,76.59) and (350.48,107.53) .. (350.48,145.7) .. controls (350.48,183.86) and (319.54,214.8) .. (281.38,214.8) .. controls (243.21,214.8) and (212.28,183.86) .. (212.28,145.7) -- cycle ;
        %Straight Lines [id:da3662707346073586] 
        \draw [color={rgb, 255:red, 0; green, 110; blue, 0 }  ,draw opacity=1 ]   (245.81,86.27) -- (170.06,145.59) ;
        %Straight Lines [id:da029924760461540578] 
        \draw [color={rgb, 255:red, 125; green, 0; blue, 0 }  ,draw opacity=1 ]   (245.81,86.27) -- (281.38,145.7) ;
        %Shape: Ellipse [id:dp04230708306582398] 
        \draw  [fill={rgb, 255:red, 0; green, 0; blue, 0 }  ,fill opacity=1 ] (165.2,145.7) .. controls (165.2,143.31) and (167.13,141.37) .. (169.52,141.37) .. controls (171.91,141.37) and (173.84,143.31) .. (173.84,145.7) .. controls (173.84,148.08) and (171.91,150.02) .. (169.52,150.02) .. controls (167.13,150.02) and (165.2,148.08) .. (165.2,145.7) -- cycle ;
        %Shape: Circle [id:dp8571155865101912] 
        \draw  [fill={rgb, 255:red, 0; green, 0; blue, 0 }  ,fill opacity=1 ] (277.06,145.7) .. controls (277.06,143.31) and (278.99,141.37) .. (281.38,141.37) .. controls (283.77,141.37) and (285.7,143.31) .. (285.7,145.7) .. controls (285.7,148.08) and (283.77,150.02) .. (281.38,150.02) .. controls (278.99,150.02) and (277.06,148.08) .. (277.06,145.7) -- cycle ;
        
        % Text Node
        \draw (187.93,106.05) node [anchor=north west][inner sep=0.75pt]    {$R$};
        % Text Node
        \draw (277.27,119.02) node [anchor=north west][inner sep=0.75pt]    {$r$};
        % Text Node
        \draw (121.39,152.16) node [anchor=north west][inner sep=0.75pt]    {$( 0,0)$};
        % Text Node
        \draw (292.42,158.27) node [anchor=north west][inner sep=0.75pt]    {$( d,0)$};
        % Text Node
        \draw (234,132.4) node [anchor=north west][inner sep=0.75pt]    {$\mathbf{A}$};
        
    \end{tikzpicture}
    \end{center}

    Here, two circles of radius $R$ and $r$ centered at $(0, 0)$ and $(d, 0)$ respectively intersect each other at two points, overlapping in a lens-shaped region $\mathbf A$. We want to know that if $R = r = 1$, what is the separation $d = \mathpzc{D}_{\text{DHA}}$ such that the area of $\mathbf A$ is $\frac\pi2$?

\subsection{Area of a Circular Segment}
    
    A circular segment is a region like the region $\mathbf S$ shown below, such that $\theta < \pi$. Note that $\theta$ in the case represents the angle $\angle AOB$. We will need a formula for its area in terms of $h$ and $R$ since it will prove useful later. 
    
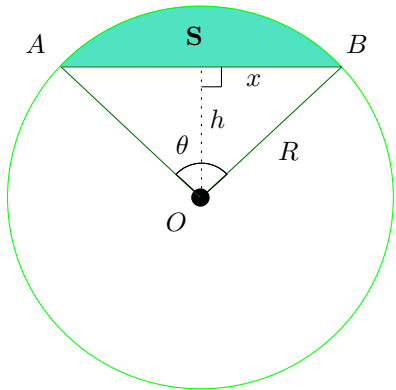
\begin{center}
    \begin{tikzpicture}[x=0.75pt,y=0.75pt,yscale=-1,xscale=1]
        %uncomment if require: \path (0,300); %set diagram left start at 0, and has height of 300
        
        %Shape: Chord [id:dp2266027266424211] 
        \draw  [color={rgb, 255:red, 80; green, 227; blue, 194 }  ,draw opacity=1 ][fill={rgb, 255:red, 80; green, 227; blue, 194 }  ,fill opacity=1 ] (99.63,79.52) .. controls (117.18,61) and (142,49.45) .. (169.52,49.45) .. controls (197.03,49.45) and (221.85,60.99) .. (239.39,79.5) -- cycle ;
        %Shape: Ellipse [id:dp3499342112834707] 
        \draw  [color={rgb, 255:red, 0; green, 255; blue, 0 }  ,draw opacity=1 ] (73.27,145.7) .. controls (73.27,92.54) and (116.36,49.45) .. (169.52,49.45) .. controls (222.68,49.45) and (265.77,92.54) .. (265.77,145.7) .. controls (265.77,198.85) and (222.68,241.94) .. (169.52,241.94) .. controls (116.36,241.94) and (73.27,198.85) .. (73.27,145.7) -- cycle ;
        %Straight Lines [id:da3662707346073586] 
        \draw [color={rgb, 255:red, 0; green, 110; blue, 0 }  ,draw opacity=1 ]   (240,80) -- (170.06,145.59) ;
        %Shape: Ellipse [id:dp04230708306582398] 
        \draw  [fill={rgb, 255:red, 0; green, 0; blue, 0 }  ,fill opacity=1 ] (165.2,145.7) .. controls (165.2,143.31) and (167.13,141.37) .. (169.52,141.37) .. controls (171.91,141.37) and (173.84,143.31) .. (173.84,145.7) .. controls (173.84,148.08) and (171.91,150.02) .. (169.52,150.02) .. controls (167.13,150.02) and (165.2,148.08) .. (165.2,145.7) -- cycle ;
        %Straight Lines [id:da23455162569900756] 
        \draw [color={rgb, 255:red, 0; green, 110; blue, 0 }  ,draw opacity=1 ]   (100,80) -- (169.52,145.7) ;
        %Straight Lines [id:da22297551619067613] 
        \draw [color={rgb, 255:red, 0; green, 110; blue, 0 }  ,draw opacity=1 ]   (100,80) -- (240,80) ;
        %Shape: Arc [id:dp47972504186801523] 
        \draw  [draw opacity=0] (157.38,134.02) .. controls (160.52,130.59) and (165.04,128.43) .. (170.06,128.43) .. controls (175.1,128.43) and (179.63,130.6) .. (182.77,134.06) -- (170.06,145.59) -- cycle ; \draw   (157.38,134.02) .. controls (160.52,130.59) and (165.04,128.43) .. (170.06,128.43) .. controls (175.1,128.43) and (179.63,130.6) .. (182.77,134.06) ;  
        %Straight Lines [id:da9026402447979207] 
        \draw  [dash pattern={on 0.84pt off 2.51pt}]  (170.06,145.59) -- (170,80) ;
        %Straight Lines [id:da8396810889503434] 
        \draw    (170,90) -- (180,90) ;
        %Straight Lines [id:da5334023985185067] 
        \draw    (180,80) -- (180,90) ;
        
        % Text Node
        \draw (207.03,116.2) node [anchor=north west][inner sep=0.75pt]    {$R$};
        % Text Node
        \draw (156,113.4) node [anchor=north west][inner sep=0.75pt]    {$\theta $};
        % Text Node
        \draw (173,100.4) node [anchor=north west][inner sep=0.75pt]    {$h$};
        % Text Node
        \draw (191,82.4) node [anchor=north west][inner sep=0.75pt]    {$x$};
        % Text Node
        \draw (161,58.4) node [anchor=north west][inner sep=0.75pt]    {$\mathbf{S}$};
        % Text Node
        \draw (81,62.4) node [anchor=north west][inner sep=0.75pt]    {$A$};
        % Text Node
        \draw (151,152.4) node [anchor=north west][inner sep=0.75pt]    {$O$};
        % Text Node
        \draw (241,62.4) node [anchor=north west][inner sep=0.75pt]    {$B$};
        
    \end{tikzpicture}
    \end{center}

    From the diagram, we can observe such that $\cos\left(\frac{\theta}{2}\right) = \frac{h}{R}$, as well as $x^2 + h^2 = R^2$.\\
    
    Rearranging the aforementioned expressions, we can solve for $\theta$ and $x$ in terms of $R$ and $h$. And since the area of $\mathbf S$ is given by subtracting the area of sector $\sector AOB$ and triangle $\triangle AOB$, we must have it such that
    \begin{align*}
        A_{\mathbf S} &= A_{\sector AOB} - A_{\triangle AOB}\\
        &= \frac{1}{2}R^2\theta - h\cdot x \\
        &= \frac{1}{2}R^2\cdot 2\arccos\left(\frac{h}{R}\right) - h\cdot \sqrt{R^2 - h^2}\\
        &= R^2\arccos\left(\frac{h}{R}\right) - h \sqrt{R^2 - h^2} \tag*{(1.2.1)}
    \end{align*}

\subsection{Area of Lens-shaped Region $\mathbf A$}
    
    Back to the main problem, we start by finding the area of this lens-shaped region $\mathbf A$. in terms of $R$, $r$, and $d$.\\

    The equations giving the two circles are 
    \begin{align*}
        x^2 + y^2 = R^2\\
        (x-d)^2 + y^2 = r^2
    \end{align*}
    
    substituting $y^2 = R^2 - x^2$ from the first equation into the second gives us the $x$ value of the two intersection points, 
    \[(x - d)^2 + (R^2 - x^2) = r^2 \iff x = d_1 = \frac{d^2 - r^2 + R^2}{2d} \tag*{(1.3.1)}\]
    
    Now, observe the diagram below
    
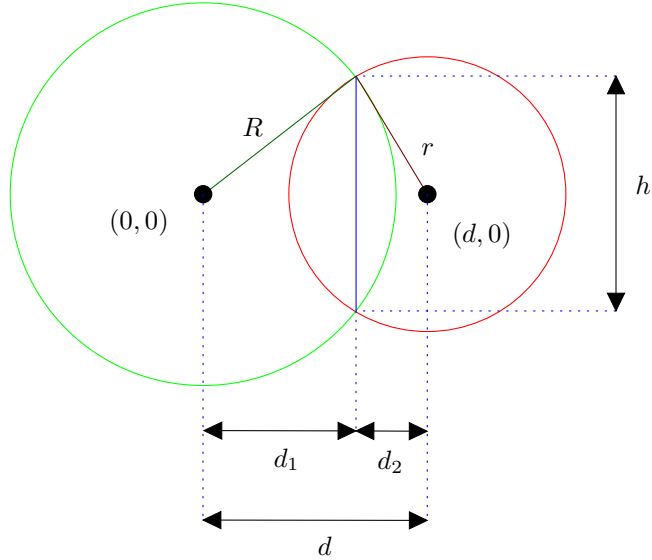
\begin{center}
    \begin{tikzpicture}[x=0.75pt,y=0.75pt,yscale=-1,xscale=1]
        %uncomment if require: \path (0,384); %set diagram left start at 0, and has height of 384
        
        %Shape: Ellipse [id:dp011685850777799711] 
        \draw  [color={rgb, 255:red, 0; green, 255; blue, 0 }  ,draw opacity=1 ] (53.27,125.7) .. controls (53.27,72.54) and (96.36,29.45) .. (149.52,29.45) .. controls (202.68,29.45) and (245.77,72.54) .. (245.77,125.7) .. controls (245.77,178.85) and (202.68,221.94) .. (149.52,221.94) .. controls (96.36,221.94) and (53.27,178.85) .. (53.27,125.7) -- cycle ;
        %Shape: Ellipse [id:dp6867065510299758] 
        \draw  [color={rgb, 255:red, 255; green, 0; blue, 0 }  ,draw opacity=1 ] (192.28,125.7) .. controls (192.28,87.53) and (223.21,56.59) .. (261.38,56.59) .. controls (299.54,56.59) and (330.48,87.53) .. (330.48,125.7) .. controls (330.48,163.86) and (299.54,194.8) .. (261.38,194.8) .. controls (223.21,194.8) and (192.28,163.86) .. (192.28,125.7) -- cycle ;
        %Straight Lines [id:da4332084795021809] 
        \draw [color={rgb, 255:red, 0; green, 110; blue, 0 }  ,draw opacity=1 ]   (225.81,66.27) -- (150.06,125.59) ;
        %Straight Lines [id:da7113031438322677] 
        \draw [color={rgb, 255:red, 125; green, 0; blue, 0 }  ,draw opacity=1 ]   (225.81,66.27) -- (261.38,125.7) ;
        %Shape: Ellipse [id:dp22273862847434578] 
        \draw  [fill={rgb, 255:red, 0; green, 0; blue, 0 }  ,fill opacity=1 ] (145.2,125.7) .. controls (145.2,123.31) and (147.13,121.37) .. (149.52,121.37) .. controls (151.91,121.37) and (153.84,123.31) .. (153.84,125.7) .. controls (153.84,128.08) and (151.91,130.02) .. (149.52,130.02) .. controls (147.13,130.02) and (145.2,128.08) .. (145.2,125.7) -- cycle ;
        %Shape: Circle [id:dp3947784311804491] 
        \draw  [fill={rgb, 255:red, 0; green, 0; blue, 0 }  ,fill opacity=1 ] (257.06,125.7) .. controls (257.06,123.31) and (258.99,121.37) .. (261.38,121.37) .. controls (263.77,121.37) and (265.7,123.31) .. (265.7,125.7) .. controls (265.7,128.08) and (263.77,130.02) .. (261.38,130.02) .. controls (258.99,130.02) and (257.06,128.08) .. (257.06,125.7) -- cycle ;
        %Straight Lines [id:da3705000121579477] 
        \draw [color={rgb, 255:red, 0; green, 0; blue, 255 }  ,draw opacity=1 ]   (225.81,66.27) -- (225.81,184.44) ;
        %Straight Lines [id:da9672658787928544] 
        \draw [color={rgb, 255:red, 0; green, 0; blue, 255 }  ,draw opacity=1 ] [dash pattern={on 0.84pt off 2.51pt}]  (149.52,125.7) -- (149.43,289.63) ;
        %Straight Lines [id:da4786117870956256] 
        \draw [color={rgb, 255:red, 0; green, 0; blue, 255 }  ,draw opacity=1 ] [dash pattern={on 0.84pt off 2.51pt}]  (225.81,184.44) -- (225.76,244.79) ;
        %Straight Lines [id:da9976965687347183] 
        \draw [color={rgb, 255:red, 0; green, 0; blue, 255 }  ,draw opacity=1 ] [dash pattern={on 0.84pt off 2.51pt}]  (261.38,125.7) -- (261.33,289.8) ;
        %Straight Lines [id:da12241437736805028] 
        \draw [color={rgb, 255:red, 0; green, 0; blue, 255 }  ,draw opacity=1 ] [dash pattern={on 0.84pt off 2.51pt}]  (355.5,66.27) -- (225.81,66.27) ;
        %Straight Lines [id:da9408699218736281] 
        \draw [color={rgb, 255:red, 0; green, 0; blue, 255 }  ,draw opacity=1 ] [dash pattern={on 0.84pt off 2.51pt}]  (355.5,184.44) -- (225.81,184.44) ;
        %Straight Lines [id:da16697705983919486] 
        \draw    (152.43,289.64) -- (258.33,289.8) ;
        \draw [shift={(261.33,289.8)}, rotate = 180.09] [fill={rgb, 255:red, 0; green, 0; blue, 0 }  ][line width=0.08]  [draw opacity=0] (8.93,-4.29) -- (0,0) -- (8.93,4.29) -- cycle    ;
        \draw [shift={(149.43,289.63)}, rotate = 0.09] [fill={rgb, 255:red, 0; green, 0; blue, 0 }  ][line width=0.08]  [draw opacity=0] (8.93,-4.29) -- (0,0) -- (8.93,4.29) -- cycle    ;
        %Straight Lines [id:da16954637531418826] 
        \draw    (152.89,244.62) -- (222.76,244.78) ;
        \draw [shift={(225.76,244.79)}, rotate = 180.13] [fill={rgb, 255:red, 0; green, 0; blue, 0 }  ][line width=0.08]  [draw opacity=0] (8.93,-4.29) -- (0,0) -- (8.93,4.29) -- cycle    ;
        \draw [shift={(149.89,244.62)}, rotate = 0.13] [fill={rgb, 255:red, 0; green, 0; blue, 0 }  ][line width=0.08]  [draw opacity=0] (8.93,-4.29) -- (0,0) -- (8.93,4.29) -- cycle    ;
        %Straight Lines [id:da4571263617839857] 
        \draw    (228.76,244.8) -- (258.29,244.95) ;
        \draw [shift={(261.29,244.96)}, rotate = 180.28] [fill={rgb, 255:red, 0; green, 0; blue, 0 }  ][line width=0.08]  [draw opacity=0] (8.93,-4.29) -- (0,0) -- (8.93,4.29) -- cycle    ;
        \draw [shift={(225.76,244.79)}, rotate = 0.28] [fill={rgb, 255:red, 0; green, 0; blue, 0 }  ][line width=0.08]  [draw opacity=0] (8.93,-4.29) -- (0,0) -- (8.93,4.29) -- cycle    ;
        %Straight Lines [id:da16009682300854666] 
        \draw    (355.5,69.27) -- (355.5,181.44) ;
        \draw [shift={(355.5,184.44)}, rotate = 270] [fill={rgb, 255:red, 0; green, 0; blue, 0 }  ][line width=0.08]  [draw opacity=0] (8.93,-4.29) -- (0,0) -- (8.93,4.29) -- cycle    ;
        \draw [shift={(355.5,66.27)}, rotate = 90] [fill={rgb, 255:red, 0; green, 0; blue, 0 }  ][line width=0.08]  [draw opacity=0] (8.93,-4.29) -- (0,0) -- (8.93,4.29) -- cycle    ;
        
        % Text Node
        \draw (167.93,86.05) node [anchor=north west][inner sep=0.75pt]    {$R$};
        % Text Node
        \draw (257.27,99.02) node [anchor=north west][inner sep=0.75pt]    {$r$};
        % Text Node
        \draw (101.39,132.16) node [anchor=north west][inner sep=0.75pt]    {$( 0,0)$};
        % Text Node
        \draw (272.42,138.27) node [anchor=north west][inner sep=0.75pt]    {$( d,0)$};
        % Text Node
        \draw (205.4,296.44) node [anchor=north west][inner sep=0.75pt]    {$d$};
        % Text Node
        \draw (183.44,252.21) node [anchor=north west][inner sep=0.75pt]    {$d_{1}$};
        % Text Node
        \draw (235.28,253.48) node [anchor=north west][inner sep=0.75pt]    {$d_{2}$};
        % Text Node
        \draw (363.91,114.87) node [anchor=north west][inner sep=0.75pt]    {$h$};
        
    \end{tikzpicture}
    \end{center}

    The area of $\mathbf A$ is the sum of two circular segments from the circles at $(0, 0)$ and $(d, 0)$, denoted as $\semicirc R$ and $\semicirc r$. They have radius and height $R$, $d_1$, and $r$, $d_2$, respectively. However, we also know from (1.3.1) that
    \[d_2 = d - d_1 = d - \frac{d^2 - r^2 + R^2}{2d} = \frac{d^2 + r^2 - R^2}{2d} \tag*{(1.3.2)}\]

    Using these results gives us the formula for the area of $\mathbf A$, 
    \begin{align*}
        A_{\mathbf A} &= A_{\semicirc \! R} + A_{\semicirc r}\\
        &= \left[ R^2 \arccos\left( \frac{d_1}{R} \right) - d_1 \sqrt{R^2 - d_1^2} \right] + \left[ r^2 \arccos\left( \frac{d_2}{r} \right) - d_2 \sqrt{r^2 - d_2^2} \right] \\
        &= \left[ R^2 \arccos\left( \frac{\frac{d^2 - r^2 + R^2}{2d}}{R} \right) - \left(\frac{d^2 - r^2 + R^2}{2d}\right) \sqrt{R^2 - \left(\frac{d^2 - r^2 + R^2}{2d}\right)^2} \right] \\
        & \phantom{==} + \left[ r^2 \arccos\left( \frac{\frac{d^2 + r^2 - R^2}{2d}}{r} \right) - \left(\frac{d^2 + r^2 - R^2}{2d}\right) \sqrt{r^2 - \left(\frac{d^2 + r^2 - R^2}{2d}\right)^2} \right]
    \end{align*}

    Simplifying gives us the formula for the area of $\mathbf A$ in terms of the radii of the two circles and their separation,
    \[A_{\mathbf A} = R^2 \arccos\left( \frac{d^2 + R^2 - r^2}{2dR} \right) + r^2 \arccos\left( \frac{d^2 - R^2 + r^2}{2dr} \right) - \frac{1}{2}\sqrt{-d^4 + 2d^2r^2 + 2d^2R^2 - r^4 + 2r^2R^2 - R^4} \tag{1.3.3} \]

\section{Overlapping Area}

\subsection{Simplifying to the Transcendental Equation}

    Now assume that for our overlapping circles, $R = r = 1$ (so they are unit circles). If both circles overlap by exactly half of each other's, that means $A_{\mathbf A} = \frac{\pi}{2}$. We can plug these into (1.3.3), and simplify over positive reals (distance cannot be negative lol) to get
    
    \begin{align*}
        \frac{\pi}{2} &= \arccos\left( \frac{d^2 + 1 - 1}{2d} \right) + \arccos\left( \frac{d^2 - 1 + 1}{2d} \right) - \frac{1}{2}\sqrt{-d^4 + 2d^2 + 2d^2 - 1 + 2 - 1}\\
        &= 2\arccos\left( \frac{d}{2} \right) - \frac{1}{2}\sqrt{4d^2 - d^4} \\
        &= 2\arccos\left( \frac{d}{2} \right) - \frac{1}{2}d\sqrt{4 - d^2} \tag{2.1.1}
    \end{align*}

    Numerically solving (2.1.1) yields $d = \mathpzc{D}_{\text{DHA}} = 0.8079455\dots$, but there is no closed form known.\cite{WolframWebsite-Circle}

\subsection{Kepler's Equation}

    We now begin by taking a look at an equation in astrophysics, called Kepler’s equation,
    $$M = E - e\sin(E)$$
    Here, $M$ represents the mean anomaly (a parameterization of time), and $E$ the eccentric anomaly (a parameterization of polar angle) of a celestial body orbiting on an elliptical path with eccentricity $e \in [0, 1)$.\\

    If we write Kepler’s equation in a less astrophysics-oriented manner, we get
    $$x = y - a \sin(y)$$
    Let's use this version from now on. For $x$ that aren't multiples of $\pi$, Kepler's equation has a unique solution for $y$ for all $a \in [0, 1)$. However, since this equation is transcendental, inverting the equation to obtain the closed form of $y$ given any arbitrary $x$ comes with great difficulty. We define the “solution” to this equation for $y$ as the "Kepler E" function $\mathbf E(a, x)$,\cite{KeplerEFunction} which is equivalent to any given series solution to the equation. In particular, let us define
    \[\mathbf{E}(a, x) = x + \sum_{n=1}^\infty \frac{2}{n} J_n(n a) \sin(n x) \tag{2.2.1}\]
    for $\mathbb{R} \ni |a| < 1$ (we don't care about astrophysics, only convergence matters) and $x\in[-\pi, \pi]$, where $J_n(x)$ is the Bessel function of the first kind.\cite{WolframWebsite-Kepler} \\

\subsection{Beta Functions}

    Now, we introduce a multitude of beta function variants. We have the normal beta function, incomplete beta function, and regularized beta function, given by, respectively, 
    \[ \operatorname{B}(a, b) = \frac{\Gamma(a)\, \Gamma(b)}{\Gamma(a+b)} = \int_0^1 t^{a-1} (1-t)^{b-1} \text{ d}t \tag{2.3.1} \]
    \[ \operatorname{B}(x; a, b) = \int_0^x t^{a-1} (1-t)^{b-1} \text{ d}t \tag{2.3.2} \]
    \[ I_x(a, b) = \frac{\operatorname{B}(x; a, b)}{\operatorname{B}(a, b)} \tag{2.3.3} \]
    The domain restriction follows after an analytical continuation of the beta function integrals, with $a, b\in \mathbb{C} \setminus 0\cup\mathbb{Z}^-$, and $x\in\mathbb{C}$. Then, we introduce the inverse regularized beta function, which serves as a quantile function for beta distributions in statistics, defined as follows. For 
    \[ z = I_x(a, b) \iff I^{-1}_z(a, b) = x \tag{2.3.4} \]
    Since it is also a quantile function, we restrict the variables $z, a, b \in \mathbb R$ to $0\le z\le 1$, and $a, b\ge 0$.

\subsection{Closed Form}
        
    Okay, how do all these new functions help? Well, the goal is to notice that
    \[
    \operatorname{B}\left( x; \frac 1 2, \frac 3 2 \right) = 
    \int_0^x \sqrt{\frac{1}{t} - 1} \text{ d}t = 
    \sqrt{x - x^2} - \arctan\left( \sqrt{\frac{1-x}{x}} \right) + \frac{\pi}{2} = 
    \sqrt{x - x^2} + \arcsin\left( \sqrt x \right) 
    \tag{2.4.1}
    \]
    where the last step follows through some trivial trig identities that can be found by drawing right triangles and phase shifting (we are simplifying over the reals).\newpage %hardcoding bruh
    
    Then, if we let $x=\sin^2(z)$, while restricting $0 \le z \le \frac{\pi}{2}$ (to simplify absolute values), we can see that 
    \[
    \operatorname{B}\left( \sin^2(z); \frac 1 2, \frac 3 2 \right) = 
    \sqrt{\sin^2(z) - \left(\sin^2(z)\right)^2} + \arcsin\left( \sqrt{\sin^2(z)} \right) = 
    \cos(z)\sin(z) + z = 
    z + \frac{1}{2}\sin(2z)
    \tag{2.4.2}
    \]
    
    Setting the entire equation equal to some dummy $y$ allows us to invert it using the Kepler E function,\cite{MSE-Dottie}
    $$
    y = 
    z + \frac{1}{2}\sin(2z) = 
    \frac{\color{red} 2z}{2} + \frac{1}{2}\sin({\color{red} 2z}) 
    \iff 
    2y = 
    {\color{red} 2z} - (-1\cdot \sin({\color{red} 2z}))
    $$
    \[ \implies z = \frac{\mathbf{E}(-1, 2y)}{2} \tag{2.4.3} \]

    However, if 
    \begin{align*}
        \operatorname{B}\left( \sin^2(z); \frac 1 2, \frac 3 2 \right) &= \operatorname{B}\left( \frac 1 2, \frac 3 2 \right) \cdot I_{\sin^2(z)}\left( \frac 1 2, \frac 3 2 \right) \\
        &= \frac{\pi}{2} I_{\sin^2(z)}\left( \frac 1 2, \frac 3 2 \right) = y \\
        \implies \sin^2(z) &= I^{-1}_{\frac{2y}{\pi}}\left( \frac{1}{2}, \frac{3}{2} \right) \\
        \implies z &= \arcsin\left(\sqrt{ I^{-1}_{\frac{2y}{\pi}}\left( \frac{1}{2}, \frac{3}{2} \right) }\right)
    \end{align*}
    this means that \[\frac{\mathbf{E}(-1, 2y)}{2} = \arcsin\left(\sqrt{ I^{-1}_{\frac{2y}{\pi}}\left( \frac{1}{2}, \frac{3}{2} \right) }\right) \tag{2.4.4} \]

    Substituting another dummy variable $x=2y$, we therefore obtain
    \[
    \mathbf{E}(-1, x) = 
    2\arcsin\left(\sqrt{ I^{-1}_{\frac{x}{\pi}}\left( \frac{1}{2}, \frac{3}{2} \right) }\right) = 
    \operatorname{archav}\left( I^{-1}_{\frac{x}{\pi}}\left( \frac{1}{2}, \frac{3}{2} \right) \right)
    \tag{2.4.5}
    \]
    $$
    \implies 
    I^{-1}_{\frac{x}{\pi}}\left( \frac{1}{2}, \frac{3}{2} \right) = 
    \sin^2\left( \frac{\mathbf{E}(-1, x)}{2} \right) = 
    \operatorname{hav}(\mathbf{E}(-1, x))
    $$
    where $\operatorname{hav}(x)$ and $\operatorname{archav}(x)$ are the haversine and inverse haversine functions, respectively.\\

    Going back to equation (2.1.1), we can make a series of substitutions with dummy variables as follows
    \begin{align*}
        \frac{\pi}{2} &= 2 \arccos\left(\frac{d}{2}\right) - \frac{d}{2}\sqrt{4-d^2} \\
        \implies \pi &= 4 \arccos\left(\alpha \right) - 4\sqrt{\alpha^2 - \alpha^4},\qquad d=2\alpha \\
        &= 2\pi - 4 \arcsin\left(\alpha \right) - 4\sqrt{\alpha^2 - \alpha^4}\\
        &= 4 \arcsin\left(\sqrt{\beta} \right) + 4\sqrt{\beta - \beta^2}, \qquad \alpha^2 = \beta\\
        \implies \frac{\pi}{4} &= \arcsin\left(\sqrt{\beta} \right) + \sqrt{\beta - \beta^2} \tag{2.4.6}
    \end{align*}
    
    Notice how our expression (2.4.6) is homologous to the beta function expression in (2.4.1). Therefore, we have
    \begin{align*}
        \frac{\pi}{4} &= \arcsin\left(\sqrt{\beta} \right) + \sqrt{\beta - \beta^2} \\
        &= \operatorname{B}\left( \beta; \frac 1 2, \frac 3 2 \right)\\
        &= \operatorname{B}\left( \sin^2(\gamma); \frac 1 2, \frac 3 2 \right), \qquad \beta = \sin^2(\gamma)\\
        &= \gamma + \frac{1}{2}\sin(2\gamma)\\
        \implies \gamma &= \frac{\mathbf{E}\left(-1, 2\cdot \frac{\pi}{4}\right)}{2} = \frac{\mathbf{E}\left(-1, \frac{\pi}{2}\right)}{2} \tag{2.4.7}
    \end{align*}

    By reversing all substitutions, we can get an expression involving the Kepler E function for $d = \mathpzc{D}_{\text{DHA}}$. Then, using (2.4.5), this allows us to conclude that $\mathpzc{D}_{\text{DHA}}$ has the following closed forms
    \begin{align*}
        \mathpzc{D}_{\text{DHA}} &= 2\sin\left(\frac{\mathbf{E}\left(-1, \frac{\pi}{2}\right)}{2}\right) \tag{2.4.8}\\
        &= 2\sin\left(\frac{1}{2}\operatorname{archav}\left( I^{-1}_{\frac{1}{2}}\left( \frac{1}{2}, \frac{3}{2} \right) \right)\right) \tag{2.4.9}\\
        &= 2\sqrt{ I^{-1}_{\frac{1}{2}}\left( \frac{1}{2}, \frac{3}{2} \right) }\tag{2.4.10}
    \end{align*}

\section{Conclusion}

\subsection{Result}

    The separation between the centers of two unit circles such that their overlapping area is exactly half of each’s area $\mathpzc{D}_{\text{DHA}}$, can be expressed as
    $$
    \mathpzc{D}_{\text{DHA}} 
    = 2\sin\left(\frac{\mathbf{E}\left(-1, \frac{\pi}{2}\right)}{2}\right)
    = 2\sin\left(\frac{1}{2}\operatorname{archav}\left( I^{-1}_{\frac{1}{2}}\left( \frac{1}{2}, \frac{3}{2} \right) \right)\right)
    = 2\sqrt{ I^{-1}_{\frac{1}{2}}\left( \frac{1}{2}, \frac{3}{2} \right) }
    $$
    where $\mathbf{E}\left(a, x\right)$ is the Kepler E function, $\operatorname{archav}(x)$ is the inverse haversine function, and $I_z^{-1}(a, b)$ is the inverse regularized beta function. Ultimately, I find (2.4.10) a particularly nice result.\\

    Numerically, this is 
    $$\mathpzc{D}_{\text{DHA}} = 0.80794550659903441863792348013263088580447192914819684450019520346774109994259070700248678\dots$$

\subsection{Acknowledgements}

    I would like to thank Tyma Gaidash for their insight in manipulating the beta function in their MSE post.
    
\printbibliography

\end{document}